\documentclass{amsart}

\RequirePackage[OT1]{fontenc}
\RequirePackage{amsthm,amsmath}
\RequirePackage[numbers]{natbib}
\RequirePackage[colorlinks,citecolor=blue,urlcolor=blue]{hyperref}

\usepackage{amsfonts}
\newtheorem{theorem}{Theorem}
\newtheorem{lemma}{Lemma}
\newtheorem{corollary}{Corollary}

\newtheorem{proposition}{Proposition}
\newtheorem{definition}{Definition}

\newcommand{\ck}{\begin{flushright} \rule{2mm}{2mm} \end{flushright}}

\newcommand{\rr}{{\mathbb R}}
\newcommand{\RR}{\mathbb{R}}

\newcommand{\DD}{\mathbb{D}}

\newcommand{\JJ}{\mathbb{J}}
\newcommand{\MM}{\mathbb{M}}
\newcommand{\SJ}{\mathcal{S}_{\mathbb{J}_1}}
\newcommand{\SM}{\mathcal{S}_{\mathbb{M}_1}}

\renewcommand{\(}{\left(}
\renewcommand{\)}{\right)}

\numberwithin{equation}{section}
\theoremstyle{plain}

\begin{document}
\title[]{Limit theorems for continuous-time random walks with continuous paths}

\author{Marcin Magdziarz  }
\address{Faculty of Pure and Applied Mathematics\\
Hugo Steinhaus Center\\
Wroclaw University of Science and Technology, Poland \\
}
\email{marcin.magdziarz@pwr.wroc.pl}

\author{Piotr \.Zebrowski}

\address{Advanced Systems Analysis (ASA) Program\\
International Institute for Applied Systems Analysis (IIASA), Laxenburg, Austria\\
}

\begin{abstract}
Continuous-time random walks are typically defined in the way that their trajectories are discontinuous step functions. This may be an unwelcome feature from the point of view of application of these processes to model certain physical and biological phenomena, where continuous paths are observed. In this article we propose alternative definition of continuous-time random walks with continuous trajectories. We also give the functional limit theorem for sequence of such random walks. In general case this result requires the use of strong Skorohod $\mathbb{M}_1$ topology instead of Skorohod $\mathbb{J}_1$ topology, which is usually used in limit theorems for ordinary continuous-time random walk processes. We also give additional conditions under which convergence of sequence of considered random walks holds in the $\mathbb{J}_1$ topology.
\end{abstract}

\maketitle

\section{Introduction}
Continuous-time random walk (CTRW) is a stochastic process determined uniquely by $\mathbb{R}^d$-valued random vectors ${Y_1}$, ${Y_2},\dots$ representing consecutive jumps of the random walker,
and $\rr_+$-valued random variables $J_1,J_2,\dots$ representing waiting times between jumps.
Every trajectory of CTRW is a step function with intervals and jumps equal to $J_i$ and
$Y_i$, respectively. Taking $J_1=J_2= {\dots}  = 1$ we obtain the classical random walk process.
CTRWs were introduced for the first time in the pioneering work of Montroll and Weiss \cite{Weiss}.
Since then they became one of the most popular and useful models in statistical physics \cite{Sokolov_book}.
Their first spectacular application can be found in \cite{Montroll}, in which
CTRW with heavy-tailed waiting times was used as a model of charge carrier transport in amorphous semiconductors.
Today CTRWs are well established mathematical models, particularly
attractive in the modeling of anomalous dynamics characterized by
nonlinear in time mean square displacement $Var(X(t))\sim t^\alpha$, $\alpha\neq 1$, see \cite{Metzler_Klafter} and references therein.

The main issues that arise in the mathematical studies of
CTRWs are limit theorems and governing equations
describing evolution in time of the corresponding probability density functions. There is an extensive literature
in this field: general results for the scaling limits of CTRWs can be found in \cite{gammagt1,StrakaHenry,coupleCTRW,CTRW1,Sikorskii}.
Governing equations for the densities of the CTRW limits and the related fractional Cauchy problems were analyzed in \cite{Baeumer,Sikorskii,Jurlewicz,Nane}.
Some recent results for particular classes of correlated and coupled CTRWs with the corresponding
Langevin picture can be found in \cite{Jurlewicz2,MM1,MM2,MM3}

The trajectories of CTRW are step functions, thus they are discontinuous.
However, the usual physical and biological requirement for a mathematical model
is to have continuous realizations. The remedy
for this problem is to apply the standard linear interpolation
to the trajectories of CTRW. As a result one obtains a linearly interpolated
CTRW with continuous, piecewise linear, trajectories -- a proper model of physical system.
Although much is known about asymptotic properties of CTRWs,
there are no such results for linearly interpolated CTRWs.
The only one exception that we are aware of, 
are  L\'evy walks which were recently studied in \cite{MPPP}.

In the case of standard random walk
its linearly interpolated version
has exactly the same limit (in the $M_1$ topology) as the random walk itself, see Corollary 6.2.1 in. \cite{Whitt}.
As our results will show, the situation can be drastically different for CTRWs and
their linearly interpolated counterparts -- the limits can differ significantly.

In this paper we intorduce a modification of CTRW. Namely, we use linear interpolation
in order to make the trajectories of CTRW continuous, which is crucial
in physical and biological applications.
We prove functional limit theorems for such CTRWs with continuous paths.
To the best of our knowledge, this is the first systematic
study of linearly interpolated CTRWs and their limits.

\section{The model}
CTRW is given by the sequence $\{(Y_k, J_k), k \geq 1\} \equiv \{(Y_k, J_k)\}$ of random vectors in $\RR^d \times [0, \infty)$, $d \geq 1$. The random vectors $Y_k \in \RR^d$ represent the successive jumps of the walker while $J_k$ are the waiting times separating moments of jumps.  We define processes of cumulative jumps and waiting times 
\begin{equation} \label{def_S_i_T}
S(t) \stackrel{df}{=} \sum_{k=1}^{[t]} Y_k \ , \qquad T(t) \stackrel{df}{=} \sum_{k=1}^{[t]} J_k \ , \quad t\geq 0
\end{equation} 
and the counting process
\begin{equation} \label{def_N} 
N(t) \stackrel{df}{=} \max\left\{k \geq 0 \ : \ T(k) \leq t \right\}, \quad t\geq 0.
\end{equation} 
Quantities $S(n)$, $T(n)$ and $N(t)$ can be interpreted respectively as the position of the walker after $n$-th jump, the moment of the $n$-th jump and the number of the last jump before time $t$.

CTRW process $X$ generated  by the sequence $\{(Y_k, J_k)\}$ is the process $S$ subordinated by counting process $N$, i.e.
\begin{equation} \label{def_CTRW}
X(t) \stackrel{df}{=} S(N(t)) = \sum_{k=1}^{N(t)} Y_k,\quad t\geq 0.
\end{equation}
The sequence $\{(Y_k, J_k)\}$ generates also the so-called overshooting CTRW process $\tilde X$ defined as
\begin{equation} \label{def_OCTRW}
\tilde X(t) \stackrel{df}{=} S(N(t)+1) = \sum_{k=1}^{N(t)+1} Y_k, \quad t\geq 0.
\end{equation}
Quantities $X(t)$ and $\tilde X(t)$ give the last position of the walker before time $t$ and the position of the walker just after first jump after time $t$, respectively. Similarly, $T(N(t))$ is a moment of the last jump before time $t$ and $T(N(t)+1)$ is the moment of the first jump after time $t$.\\

From the definition of CTRW process it follows that its trajectories are the \textit{c\'adl\'ag} (i.e. left-continuous and having right limits) step functions. However sequence $\{(Y_k, J_k)\}$ generates also a process $\bar X$ with continuous trajectories
\begin{equation} \label{def_CPCTRW}
\bar X(t) \stackrel{df}{=} S(N(t)) + \frac{S(N(t)+1) - S(N(t))}{T(N(t)+1) - T(N(t))}\(t - T(N(t))\) \quad t\geq 0.
\end{equation}
Indeed, the trajectories of process $\bar X$ are piecewise linear. We call $\bar X$ \textbf{continuous time random walk with continuous paths} (or CPCTRW i.e. continuous-path CTRW). Such processes are the main interest of this article.\\

Observe that for each fixed $\omega \in \Omega$ the trajectory of $\bar X(\omega, \cdot)$ arises from the trajectory of $X(\omega, \cdot)$ by replacing each jump of $X(\omega, \cdot)$, say at time $t$ ,  with the segment linking the 
points $(S(N(t)),T(N(t)))$ and $(S(N(t)+1),T(N(t))+1)$. 
Thus we may write  $\bar X(\omega, t) = f(X(\omega, \cdot))(t)$, where $f$ is some appropriate mapping. In the following section we will define this mapping $f$ and investigate its properties.

\section{The mapping modifying trajectories. Definition and basic properties.}

Let $\DD([0,\infty), \RR^d)$, $d \geq 1$, denote the space of \textit{c\'adl\'ag} functions defined on $[0,\infty)$ and taking values in $\RR^d$. In this section we define mapping $f \ : \ \DD([0, \infty), \RR^d) \mapsto \DD([0, \infty), \RR^d)$ which modifies any \textit{c\'adl\'ag} function by replacing each "stair" of this function with the segment (in the way as it was for $X$ and $\bar{X}$ above). 
It is important to stress that we say that function $x\in\DD([0,\infty), \RR^d)$ "has a stair", say at time $t$, only if it is constant on some open interval $(t-\epsilon,t)$ and $t$ is a discontinuity point (moment of jump) of $x$.
The mapping $f$ is designed to modify only the "stairs" of the underlying function $x$, thus $f(x)$ may still have discontinuities. Point $\tau$ will be the discontinuity point of $f(x)$ whenever $\tau \in disc(x)$ and $x$ is not constant on the interval $[\tau - \delta, \tau)$ for any $\delta > 0$.

In order to give a convenient definition of the mapping $f$, we need to be able to check if for any fixed $t\geq 0$ the function $x$ is constant in some neighbourhood of $t$. We also need to identify the ends of the interval on which $x$ is constant.

\begin{definition}
For any $x \in \DD([0, \infty), \RR^d)$ let
\begin{eqnarray}	
\eta_x(t) &\stackrel{def}{=}& \sup\{s<t \ : \ x(s) \neq x(t)\}, \quad t\geq 0, \label{def_eta} \\  
\theta_x(t) &\stackrel{def}{=}& \inf\{s \geq t \ : \ x(s) \neq x(t-)\}, \quad t\geq 0.    \label{def_theta}
\end{eqnarray}
\end{definition}

One easily sees that if $x$ is constant on some interval, say $(a, b)$, such that $a < t < b$ and this interval is the longest possible (i.e. $x$ is not constant on any interval $(a', b')$ where $a' < a$, $b' > b$) then $\eta_x(t) = a$ and $\theta_x(t) = b$. Further on the following lemma will prove useful.

\begin{lemma} \label{lem_eta_theta}
$\eta_x(t)$ and $\theta_x(t)$ have the following properties:
\begin{itemize}
	\item[(i)] $\eta_x(t) \leq t$ and $\eta_x(t) < t$ if and only if $x$ is constant on interval $[\eta_x(t), t]$;
	\item[(ii)] $\theta_x(t) \geq t$ and $\theta_x(t) > t$ if and only if $t \notin disc(x)$ and $x$ is constant on interval $[t, \theta_x(t))$;
\end{itemize}
\end{lemma}

\noindent \textbf{Proof.} 
From (\ref{def_eta}) it immediately follows that $\eta_x(t) \leq t$ and $\eta_x(t) < t$ if and only if there exist $\delta > 0$ such that $x$ is constant on interval $[t-\delta, t]$. One easily checks that if $\delta_0$ is the biggest of all such $\delta > 0$, then $t-\delta_0 = \eta_x(t)$ and $(i)$ holds. Similarly, from (\ref{def_theta}) we have that $\theta_x(t) \geq t$ and $\theta_x(t) > t$ if and only if $t \notin disc(x)$ and $x$ is constant on interval $[t, t+\delta)$ for some $\delta > 0$. Again we take $\delta_0$ being the biggest $\delta >0$ for which $x$ is constant on interval $[t, t+\delta)$ and check that $\theta_x(t) = t + \delta_0$. 
\ck

Having identified possible ends of "stairs" we define the mapping $f$.

\begin{definition}
The mapping $f : \DD([0, \infty), \RR^d) \mapsto \DD([0, \infty), \RR^d)$ is defined as
\begin{equation} \label {def_f}
f(x)(t) \stackrel{def}{=} 
\left\{ 
\begin{array}{ll} 
x\(\eta_x(t)\) + \frac{x\(\theta_x(t)\) - x\(\eta_x(t)\)}{\theta_x(t) - \eta_x(t)}\(t - \eta_x(t)\) & , \ \theta_x(t) > \eta_x(t) \\ 
x(t) & , \ \theta_x(t) = \eta_x(t) = t
\end{array} 
\right. .
\end{equation}
where $x \in \DD([0, \infty), \RR^d)$ and $t \geq 0$.
\end{definition}

The next lemma proves that mapping $f$ indeed modifies only the "stairs" of function $x$.

\begin{lemma} \label{lem_f}
$f(x)(t) \neq x(t)$ if and only if $x$ is constant on the interval $[\eta_x(t), \theta_x(t))$ and $\theta_x(t) \in disc(x)$.
\end{lemma}

\noindent \textbf{Proof.} 
By (\ref{def_f}) we have that $f(x)(t) \neq x(t)$ whenever (a) $x(\theta_x(t)) \neq x(\eta_x(t))$ and (b) $t > \eta_x(t)$. From Lemma \ref{lem_eta_theta} $(i)$ and (b) it follows that $x$ is constant on $[\eta_x(t), t]$ and $x(t) = x(\eta_x(t))$. Then we have also that $t < \theta_x(t)$. Indeed, if $t = \theta_x(t)$ then $x(\theta_x(t)) = x(t) = x(\eta_x(t))$ which contradicts (a). Thus by Lemma \ref{lem_eta_theta} $(ii)$ function $x$ is constant on interval $[\eta_x(t), \theta_x(t))$ where $\eta_x(t) < t < \theta_x(t)$. Then $x(\theta_x(t)-) = x(\eta_x(t))$, so (a) implies that $\theta_x(t) \in disc(x)$.

Now assume that $x$ is constant on interval $[\eta_x(t), \theta_x(t))$ 
and $\theta_x(t) \in disc(x)$. Then $x$ is constant on $[\eta_x(t), t]$ and $x(\eta_x(t)) = x(t)$. Moreover, since $x$ is constant on $[t, \theta_x(t))$ and $\theta_x(t) \in disc(x)$ it follows that $x(\theta_x(t)) \neq x(t)$. Hence $x(\theta_x(t)) - x(\eta_x(t)) \neq 0$ and $f(x)(t) \neq x(t)$. \ck

We conclude this section showing that mapping $f$ is measurable, i.e. that
$$
\forall A \in \SM \quad \{x : f(x) \in A\} \in \SJ,
$$  
where $\sigma$-fields $\SJ$ and $\SM$ are the Borel $\sigma$-fields generated by the open sets in the Skorohod $\JJ_1$ and $S\MM_1$ topologies, respectively (see \cite{Whitt}). We will use the proposition below
\begin{proposition}[Th. 11.5.2. \cite{Whitt}] \label{prop_s_f}
Borel $\sigma$-fields $\SJ$ and $\SM$ coincides with the Kolmogorov $\sigma$-field generated by the coordinate projections.
\end{proposition}

\begin{lemma} \label{lem_measurab}
The following mappings are measurable:
\begin{itemize}
	\item[(1)] $x(\cdot) \mapsto \eta_x(\cdot), \quad  x\in \DD([0, \infty), \RR^d)$;
	\item[(2)] $x(\cdot) \mapsto \theta_x(\cdot), \quad  x\in \DD([0, \infty), \RR^d)$;
	\item[(3)] $x(\cdot) \mapsto x(\eta_x(\cdot)), \quad  x\in \DD([0, \infty), \RR^d)$;
	\item[(4)] $x(\cdot) \mapsto x(\theta_x(\cdot)), \quad  x\in \DD([0, \infty), \RR^d).$
\end{itemize}
\end{lemma}

\noindent \textbf{Proof.} 
We begin with showing that the  sets $\{x : x=const \mbox{ on } [a,b]\}$, $0 \leq a<b$, and $\{x : t \in disc(x)\}$ for fixed $t \geq 0$ are measurable. 

Fix arbitrary $0 \leq a < b$ and for any $n$ let $Q(n) = \{ m/n: \ m\in \mathbb{Z}\}$. Then
$$\{x : x=const \mbox{ on } [a,b]\}$$
$$ = \bigcap_{n=1}^\infty \bigcup_{q \in Q(n)} \bigcap_{k=1}^\infty \left( \left\{x: x(t_i) \in \left(\frac{nq-1}{n}, \ q \right), \ t_i = a + \frac{b-a}{k}i, \ i=0, \ldots, k \right\}\right),$$ 
thus it is measurable. Now fix arbitrary $t \geq 0$. 
$$\{x : t \in disc(x)\}$$
$$ = \bigcup_{n=1}^\infty \bigcup_{q \in Q(n)} \bigcap_{k=1}^\infty \left( \left\{x: x(t) \geq q, \  x(t_i) < q, \ t_i = t - \frac{k-i}{nk}, \ i=0, \ldots, k-1 \right\} \right.$$
$$ \cup \left. \left\{x: x(t) \leq q, \  x(t_i) > q, \ t_i = t - \frac{k-i}{nk}, \ i=0, \ldots, k-1 \right\}\right),$$ 
so $\{x : t \in disc(x)\}$ is also measurable.

Now we show that mapping $x(\cdot) \mapsto \eta_x(\cdot)$ is measurable. By proposition \ref{prop_s_f} it is sufficient to show that the set $\{x : \eta_x(t) \leq a\}$ is measurable for any $t \geq 0$ and any $a \geq 0$. By lemma \ref{lem_eta_theta} $(i)$
$$\{x : \eta_x(t) \leq a\} = \{x : \sup\{s < t: x(s) \neq x(t)\} \leq a\}$$
$$
= \left\{ 
\begin{array}{ll} 
\{x : x=const \mbox{ on } [a,t]\} & , \ a < t \\ 
\{x : x \mbox{ is not constant on any interval } (s,t), \ s<t\} & , \ a \geq t
\end{array} 
\right. .
$$
Thus the set $\{x : \eta_x(t) \leq a\}$ is measurable since the sets $\{x : x=const \mbox{ on } [a,t]\}$ and 
$$\{x : x \mbox{ is not constant on any interval } (s,t), \ s<t\}$$
$$ = \bigcap_{k=1}^\infty \{x : x = const \mbox{ on } (t-k^{-1}, t)\}^{c}$$
are measurable.

By lemma \ref{lem_eta_theta} $(ii)$ it follows also that the set $\{x : \theta_x(t) \geq a\}$ is measurable since 
$$\{x : \theta_x(t) \geq a\} = \{x : \inf\{s > t: x(s) \neq x(t-)\} \geq a\}$$
$$
= \left\{ 
\begin{array}{ll} 
\DD & , \ a < t \\ 
\{x : t \notin disc(x)\} \cap \{x : x=const \mbox{ on } [t,a]\} & , \ a \geq t
\end{array} 
\right. .
$$
Finally, mappings $x(\cdot) \mapsto x(\eta_x(\cdot))$ and $x(\cdot) \mapsto x(\theta_x(\cdot))$ are measurable since all the functions from space $\DD$ are measurable and the composition of measurable mappings is again measurable. 
\ck

As an immediate consequence of Lemma \ref{lem_measurab} and the fact that addition and multiplication are measurable we get the following corollary:

\begin{corollary}
The mapping $f$ is measurable.
\end{corollary}

\section{Continuity of the mapping $f$}

\subsection{Convergence in the Skorohod $S\MM_1$ topology}
\ \\

In the next paragraph we give conditions, under which the mapping $f$ preserves convergence in the Skorohod $S\MM_1$ topology. Before doing so we introduce the notion of convergence of sequence of functions $x_n$ to some function $x$ in the $S\MM_1$ topology, which is based on the concept of convergence of the completed graphs of functions $x_n$ to the completed graph of $x$. Detailed definitions and extensive summary of properties of $(\DD([0,\infty), \RR^d), S\MM_1)$ space can be found in \cite{Whitt}. Below we outline only some necessary definitions and facts on convergence in the strong $\MM_1$ topology.\\

In $\RR^d$ we use maximum norm $\|\cdot\|$, i.e.
$$
\|a\| = \max_{1 \leq i \leq d}|a^i|, \quad a = (a^1, \ldots a^d) \in \RR^d.
$$
For $a, b \in \RR^d$ we define the standard segment
$$
[a, b] \stackrel{def}{=} \{\alpha a + (1-\alpha)b \ : \ \alpha \in [0,1]\}.
$$
For any function $x \in \DD([0, T], \RR^d)$, $T > 0$, we define its completed graph 
$$
\Gamma_x \stackrel{def}{=} \{ (z, t) \in \RR^d \times [0, T] \ : \ z \in [x(t-), x(t)] \}
$$
where $x(t-)$ denotes the left-hand side limit of $x$ in point $t$.

On any completed graph $\Gamma_x$ we define relation $"\leq"$ inducing total order on $\Gamma_x$. For $(z_1, t_1), (z_2, t_2) \in \Gamma_x$ we say that $(z_1, t_1) \leq (z_2, t_2)$ if one of the following holds
\begin{itemize}
	\item[(i)] $t_1 < t_2$ or
	\item[(ii)] $t_1 = t_2$ and $\| x(t_1-) - z_1\| \leq \| x(t_1-) - z_2\|$.
\end{itemize}

Let $A = \{(z_i, t_i) \in \Gamma_x, i=0, \ldots m\}$ be a $m+1$ point subset of $\Gamma_x$ such that
$$
(x(0), 0) = (z_0, t_0) \leq (z_1, t_1) \leq \ldots \leq (z_m, t_m) = (x(T), T).
$$
Then $A$ is called \emph{the ordered subset} of ordered completed graph $(\Gamma_x, \leq)$. We define \emph{the order-consistent distance between ordered set $A$ and $\Gamma_x$} as
$$
\hat d(A, \Gamma_x) \stackrel{def}{=} \max_{0\leq i < m} \sup\{\|(z, t) - (z_i, t_i)\| \vee \|(z, t) - (z_{i+1}, t_{i+1})\| \ : (z_i, t_i) \leq (z, t) < (z_{i+1}, t_{i+1}) \},
$$ 
where $(z_i, t_i) \in A$ and $(z, t) \in \Gamma_x$. We use $\vee$ as a symbol of maximum, i.e. $a \vee b \equiv \max\{a, b\}$.\\

Ordered subset of completed graph and the order-consistent distance are useful in characterizing the convergence in $S\MM_1$ topology.

\begin{proposition}[see th. 12.5.1. (vi) \cite{Whitt}] \label{prop_conv_M_1}
Let $x_n, x \in \DD([0, T], \RR^d)$. The sequence $x_n$ converges to $x$ in $S\MM_1$ topology if and only if for any $\varepsilon > 0$ there exist an integer $m$, an ordered subset $A = \{(z_i, t_i), i = 0, 1, \ldots, m \} \subset \Gamma_x$ such that $\hat d(A, \Gamma_x) < \varepsilon$, an integer $n_1$ and ordered subsets $A_n = \{(z_{n,i}, t_{n,i}), i = 0, 1, \ldots, m \} \subset \Gamma_{x_n}$ such that for all $n \geq n_1$ $\hat d(A_n, \Gamma_{x_n}) < \varepsilon$ and $d^*(A, A_n) < \varepsilon$, where
$$
d^*(A, A_n) \stackrel{def}{=} \max_{1 \leq i \leq m} \{\|(z_i, t_i) - (z_{n,i}, t_{n,i})\| \ : \ (z_i, t_i) \in A, (z_{n,i}, t_{n,i})\in A_n \}.
$$
\end{proposition}

The above proposition is also useful in showing the convergence of elements of space $(\DD([0, \infty), \RR^d), S\MM_1)$.
Let $x_n, x \in \DD([0, \infty), \RR^d)$ and let $x_n|_{[0, T]}, x|_{[0, T]}\in \DD([0, T], \RR^d)$ denote their restrictions to interval $[0, T]$. By theorem 12.9.3. \cite{Whitt} we have that $x_n \to x$ in the $S\MM_1$ topology if and only if $x_n|_{[0, T]} \to x|_{[0, T]}$ in the $S\MM_1$ topology for all $T \notin disc(x)$.  \\

\subsection{Continuity of mapping $f$ in the $S\MM_1$ topology}
\ \\

Now we discuss continuity properties of mapping $f$. In general, this mapping is not continuous in the $S\MM_1$ topology as the following example shows.\\

\noindent \textbf{Example 1:} Consider elements of $\DD([0, \infty), \RR)$
\begin{displaymath}
x_n(t) = \left\{ 
\begin{array}{ll} 
1  & , \ t \in [1, 2 - 1/n)\\
1 + 1/n &  , \ t \in [2 - 1/n, 2)\\
t & , \ t \in \RR \setminus [1, 2) 
\end{array} 
\right. 
\ , \quad  
x(t) = \left\{ 
\begin{array}{ll} 
1 & , \ t \in [1, 2)\\ 
t & , \ t \notin [1, 2) 
\end{array} 
\ .
\right. 
\end{displaymath}
Then $x_n \to x$ as $n \to \infty$ in the $S\MM_1$ topology (and even in the uniform convergence topology), but we have that $f(x_n) \to x$ in the $S\MM_1$ topology while $f(x) \equiv e$, where $e$ is the identity function. \\

\begin{theorem} \label{conv_f}
Let $x, x_n \in \DD([0, \infty), \RR^d)$. Assume that sequence $x_n$ converges to $x$ in the Skorohod $\JJ_1$ topology and satisfies condition that whenever for some $t \geq 0$ we have that $\eta_x(t) \leq t < \theta_x(t)$ and $\theta_x(t) \in disc(x)$ (i.e. $t$ belongs to some ``stair'' of $x$), then there exists sequence $t_n \to t$ such that 
\begin{itemize}
	\item[(i)] $\eta_{x_n}(t_n) \to \eta_x(t)$,
	\item[(ii)] $\theta_{x_n}(t_n) \to \theta_x(t)$,
	\item[(iii)] $\(x_n(\theta_{x_n}(t_n)) - x_n(\theta_{x_n}(t_n)-)\) \to \(x(\theta_x(t)) - x(\theta_x(t)-)\)$.
\end{itemize}
Then we also have the convergence $f(x_n) \to f(x)$ in the $S\MM_1$ topology.
\end{theorem}

\noindent \textbf{Remark 1.} The key feature of convergence of $x_n$ to $x$ in $\JJ_1$ topology is the following. If $x$ has a jump at $t$, then for all $n$ sufficiently large $x_n$ must have ``matching jump'' at some $t_n$ such that $t_n \to t$ and $\(x_n(t_n) - x_n(t_n -)\) \to \(x(t) - x(t-)\)$. However, this mode of convergence does not guarantee that $f(x_n)$ converges to $f(x)$ (see Example 1.). The latter convergence to hold requires additionally that whenever $x$ has a ``stair'' then for all sufficiently large $n$ functions $x_n$ must have a ``matching stair''. This is ensured by conditions $(i) - (iii)$.

\noindent \textbf{Remark 2.} In Example 1 ``stairs'' $x_n$ does not converge to the ``matching stair'' of $x$. Note that for all $t \in (1, 2)$ we can easily find sequence $t_n \to t$ such that $\eta_{x_n}(t_n) = \eta_x(t) = 1$, $\theta_{x_n}(t_n) = 2 -1/n \to \theta_x(t) = 2$ (thus points $(i)$ and $(ii)$ are satisfied), but none of such sequences satisfies assumption $(iii)$ since 
$$\|x_n(\theta_{x_n}(t_n)) - x_n(\theta_{x_n}(t_n)-)\| = 1/n \not\to \|x(\theta_x(t)) - x(\theta_x(t)-)\| = 1.$$

\noindent \textbf{Remark 3.} Observe that assumptions $\eta_x(t) < t < \theta_x(t)$ and $t_n \to t$ together with $(i)$ and $(ii)$ imply that for sufficiently large $n$ we have that $\eta_{x_n}(t_n) < t_n < \theta_{x_n}(t_n)$. By Lemma \ref{lem_eta_theta} $x_n$ is then constant on the interval $\(\eta_{x_n}(t_n), \theta_{x_n}(t_n)\)$ and the $\JJ_1$ convergence implies also that $x_n(t_n) \to x(t)$. Moreover $(iii)$ and assumption that $\theta_x(t) \in disc(x)$ yields that $\theta_{x_n}(t_n) \in disc(x_n)$ for sufficiently large $n$.\\

\noindent \textbf{Proof of Theorem \ref{conv_f}.}
Let $y_n \equiv f(x_n)$, $y \equiv f(x)$. By Th. 12.9.3 \cite{Whitt} to show convergence $y_n \to y$ in the space $(\DD([0, \infty), \RR^d), S\MM_1)$ it is sufficient to show that $y_n|_{[0, T]} \to y|_{[0, T]}$ in the relative $S\MM_1$ topology on $\DD([0, T], \RR^d)$ for all $T \notin disc(y)$.

Fix arbitrary $T \notin disc(y)$. Without loss of generality we may assume that $y(T) = x(T)$. Observe that by (\ref{def_f}) we have that $y(\theta_x(t)) = x(\theta_x(t))$ for any $t \geq 0$. Thus if $y(T) \neq x(T)$, then we show convergence $y_n|_{[0, \theta_x(T)]} \to y|_{[0, \theta_x(T)]}$ in the $S\MM_1$ topology which also implies convergence of $y_n|_{[0, T]}$ to $y|_{[0, T]}$.

Now fix arbitrary $\varepsilon > 0$. By Proposition \ref{prop_conv_M_1} we have that there exists an integer $m$ and an ordered subset $A = \{(z_i, t_i), i = 0, 1, \ldots, m \} \subset \Gamma_x$ such that $\hat d(A, \Gamma_x) < \varepsilon$. On the interval $[0, T]$ function $x$ has finitely many, say $K = K(\varepsilon)$, jumps of magnitude greater or equal to $\varepsilon$. Let
$$
G = G(\varepsilon) \stackrel{def}{=} \{\tau \in disc(x) : \|x(\tau) - x(\tau-)\| \geq \varepsilon \} = \{\tau^{(1)} <  \ldots < \tau^{(K)}\}. 
$$ 
Without loss of generality we may assume that all points
$$
\(x(\tau-), \sup\{s < \tau : x(s) \neq x(\tau -)\} \) \mbox{ and } \( x(\tau), \tau\), \qquad \tau \in G
$$
belong to the set $A$. Let 
$$
1\leq i_1 < i_2 < \ldots < i_K < m
$$ 
be the indices of points $(z_{i_l}, t_{i_l}) \in A$ such that
$$
(z_{i_l}, t_{i_l}) = \(x(\tau^{(l)}-), \sup\{s < \tau^{(l)} : x(s) \neq x(\tau^{(l)} -)\} \), \quad l = 1, 2, \ldots , K.
$$
Similarily, let 
$$
1 < j_1 < j_2 < \ldots < j_K \leq m
$$ 
be the indices of points $(z_{j_l}, t_{j_l}) \in A$ such that
$$
(z_{j_l}, t_{j_l}) = \(x(\tau^{(l)}), \tau^{(l)}\), \quad l = 1, 2, \ldots , K.
$$
Obviously we have that $i_1 < j_1 \leq i_2 < j_2 \leq i_3 < \ldots \leq i_K < j_K$.\\

Convergence of sequence $x_n$ to $x$ in $\JJ_1$ topology together with assumptions $(i) - (iii)$ implies that we can choose sequence of sets 
$$
G_n = G_n(\varepsilon) \stackrel{def}{=} \{\tau_n^{(1)} <  \ldots < \tau_n^{(K)}\} \subset [0, T] 
$$ 
such that $\forall \ 1 \leq l \leq K$
\begin{eqnarray*}
&& \left|\tau_n^{(l)} - \tau^{(l)} \right| \to 0, \\
&& \left|\sup\{s < \tau_n^{(l)} : x_n(s) \neq x_n(\tau_n^{(l)}-)\} - \sup\{s < \tau^{(l)} : x(s) \neq x(\tau^{(l)}-)\} \right| \to 0, \\
&& \left|\left(x_n(\tau_n^{(l)}) - x_n(\tau_n^{(l)}-) \right) - \left(x(\tau^{(l)}) - x(\tau^{(l)}-) \right)  \right| \to 0 
\end{eqnarray*}
as $n \to \infty$. Therefore we may choose $n_1 = n_1(\varepsilon)$ such that $\forall n \geq n_1$
\begin{eqnarray*}
&& G_n \subset disc(x_n), \\
&& \max_{1\leq l \leq K} \left|\tau_n^{(l)} - \tau^{(l)} \right| < \varepsilon, \\
&& \max_{1\leq l \leq K} \left|\sup\{s < \tau_n^{(l)} : x_n(s) \neq x_n(\tau_n^{(l)}-)\} - \sup\{s < \tau^{(l)} : x(s) \neq x(\tau^{(l)}-)\} \right| < \varepsilon, \\
&& \max_{1\leq l \leq K} \left| \left(x_n(\tau_n^{(l)}) - x_n(\tau_n^{(l)}-) \right) - \left(x(\tau^{(l)}) - x(\tau^{(l)}-) \right)  \right| < \varepsilon, \\
&& \sup_{\tau_n \in disc(x_n) \setminus G_n} \left| x_n(\tau_n) - x_n(\tau_n - ) \right| < \varepsilon. 
\end{eqnarray*}
Moreover, by Proposition \ref{prop_conv_M_1} there exists an integer $n_2$ and ordered subsets $A_n = \{(z_{n,i}, t_{n,i}), i = 0, 1, \ldots, m \} \subset \Gamma_{x_n}$ such that $\hat d(A_n, \Gamma_{x_n}) < \varepsilon$ and $d^*(A, A_n) < \varepsilon$ for all $n \geq n_2$. Let $n_0 = \max\{n_1, n_2\}$. Without loss of generality we may assume that for all $n \geq n_0$ points
$$
\(x_n(\tau_n-), \sup\{s < \tau_n : x_n(s) \neq x_n(\tau_n -)\} \) \mbox{ and } \( x_n(\tau_n), \tau_n\), \qquad \tau_n \in G_n
$$
belong to the set $A_n$ and that for indices $i_l$, $j_l$, $l \leq K,$ we have that
\begin{eqnarray*}
(z_{n, i_l}, t_{n, i_l}) &=&  \(x_n(\tau_n^{(l)}-), \sup\{s < \tau_n^{(l)} : x_n(s) \neq x_n(\tau_n^{(l)} -)\} \), \\
(z_{n, j_l}, t_{n, j_l}) &=&  \(x_n(\tau_n^{(l)}), \tau_n^{(l)} \).
\end{eqnarray*}

Now we construct the ordered subsets 
\begin{eqnarray*}
B &=& \{(z'_i, t'_i), i=0, 1, \ldots, m\} \subset \Gamma_y \\
B_n &=& \{(z'_{n,i}, t'_{n,i}), i=0, 1, \ldots, m\} \subset \Gamma_{y_n}
\end{eqnarray*}
such that 
$$ \hat d(B, \Gamma_y) < 3\varepsilon $$
and for all $n \geq n_0$
$$ \hat d(B_n, \Gamma_{y_n}) < 3\varepsilon \ ,$$
$$  d^*(B, B_n) < 3\varepsilon.$$\\

We begin with construction of the set $B$. For any $l \leq K$ let $I_l = \{i_l, \ldots, j_l\}$.

First we construct points $(z'_k, t'_k)$ for $k \in I_l$, $l = 1, \ldots, K$. Namely we take
\begin{eqnarray*}
z'_k &=& z_{i_l} + \(\frac{z_{j_l} - z_{i_l}}{j_l - i_l}\)\(k - i_l\), \\
t'_k &=& t_{i_l} + \(\frac{t_{j_l} - t_{i_l}}{j_l - i_l}\)\(k - i_l\).
\end{eqnarray*}
In other words, points $(z'_k, t'_k)$, $k \in I_l$ for some $l \in \{1, \ldots, K\}$ are the $(j_l - i_l + 1)$ points laying uniformly on the segment $[(z_{i_l}, t_{i_l}), (z_{j_l}, t_{j_l})]$.

We choose points $(z'_k, t'_k)$ with indices 
$$k \in \{0, \ldots, m\} \setminus \bigcup_{1 \leq l \leq K} I_l$$
in the way that $t'_k = t_k$ and $z'_k$ is such that $(z'_k, t'_k) = (z'_k, t_k) \in \Gamma_y$.\\

Now we show that 
$$
\hat d(B, \Gamma_y) = \max_{0\leq i < m} \sup\{\|(z', t') - (z'_i, t'_i)\| \vee \|(z', t') - (z'_{i+1}, t'_{i+1})\|\} < 3\varepsilon,
$$ 
where the supremum is taken over all $(z'_i, t'_i) \leq (z', t') < (z'_{i+1}, t'_{i+1})$. In order to do so we need to consider the following four cases:

$(1^\circ)$ Let $k, k+1 \in I_l$ for some $l \in \{1, \ldots, m\}$. Then $(z', t') \in [(z'_k, t'_k), (z'_{k+1}, t'_{k+1}))$ so there
exists $\alpha \in (0, 1]$ such that $(z', t') = \alpha (z'_k, t'_k) + (1-\alpha)(z'_{k+1}, t'_{k+1})$. Hence
$$\|(z', t') -  (z'_k, t'_k)\| = 
\|\alpha (z'_k, t'_k) + (1-\alpha)(z'_{k+1}, t'_{k+1}) - (z'_k, t'_k)\|$$
$$ = (1-\alpha) \|(z'_k, t'_k) - (z'_{k+1}, t'_{k+1})\| 
< \| (z'_{k+1}, t'_{k+1}) - (z'_k, t'_k)\| 
= \|z'_{k+1} - z'_k\| \vee |t'_{k+1} - t'_k|$$
$$ = \left\|\(z_{k+1} + \(\frac{z_{j_l} - z_{i_l}}{j_l - i_l}\)\(k + 1- i_l\)\) - \(z_k + \(\frac{z_{j_l} - z_{i_l}}{j_l - i_l}\)\(k - i_l\)\) \right\|$$
$$ \vee \left|\(t_{k+1} + \(\frac{t_{j_l} - t_{i_l}}{j_l - i_l}\)\(k + 1- i_l\)\) - \(t_k + \(\frac{t_{j_l} - t_{i_l}}{j_l - i_l}\)\(k - i_l\)\) \right|$$
$$= \left\|z_{k+1} - z_k + \frac{z_{j_l} - z_{i_l}}{j_l - i_l}\right\| \vee \left|t_{k+1} - t_k + \frac{t_{j_l} - t_{i_l}}{j_l - i_l}\right|.$$
Since $\hat d(A, \Gamma_x) < \varepsilon$ it easily follows that $\|z_{k+1} - z_k\| \leq \varepsilon$ and $|t_{k+1} - t_k| \leq \varepsilon$. Applying the same argument we have that 
$$\left\| \frac{z_{j_l} - z_{i_l}}{j_l - i_l} \right\| 
= \frac{1}{j_l - i_l} \left\| \sum_{s=i_l}^{j_l-1} (z_{s+1} - z_s) \right\|
\leq \frac{j_l - i_l - 1}{j_l - i_l} \varepsilon < \varepsilon,$$
$$\left| \frac{t_{j_l} - t_{i_l}}{j_l - i_l} \right| 
= \frac{1}{j_l - i_l} \left| \sum_{s=i_l}^{j_l-1} (t_{s+1} - t_s) \right|
\leq \frac{j_l - i_l - 1}{j_l - i_l} \varepsilon < \varepsilon.$$
Hence
$$\|(z', t') -  (z'_k, t'_k)\| 
\leq \left\|z_{k+1} - z_k + \frac{z_{j_l} - z_{i_l}}{j_l - i_l}\right\| \vee \left|t_{k+1} - t_k + \frac{t_{j_l} - t_{i_l}}{j_l - i_l}\right| < 2\varepsilon.$$
In the same way we show that $\|(z'_{k+1}, t'_{k+1}) - (z', t')\| < 2\varepsilon$.  

$(2^\circ)$ Let $k = j_l$ for some $l \in \{1, \ldots, m\}$. Observe that if $j_l = i_{l+1}$, then $k, k+1 \in I_{l+1}$ and we are in case $(1^\circ)$. Therefore we may assume that $j_l < i_{l+1}$. Then for any $(z'_k, t'_k) \leq (z', t') \leq (z'_{k+1}, t'_{k+1})$ we have that
$$\|(z', t') -  (z'_k, t'_k)\| 
\leq \|(z', t') -  (z, t)\| + \|(z, t) -  (z_k, t_k)\| + \|(z_k, t_k) -  (z'_k, t'_k)\|.$$
Note that
$$\|(z', t') -  (z, t)\| = \|(f(x)(t), t) - (x(t), t)\| = \|f(x)(t) - x(t)\|$$
$$ \leq \| x(\theta_x(t)) - x(\theta_x(t) - ) \| < \varepsilon.$$
The last inequality is the consequence of the assumption that $k = j_l$ for some $l$. Indeed, if $\| x(\theta_x(t)) - x(\theta_x(t) - ) \| \geq \varepsilon$ then $\theta_x(t) = \tau^{(l+1)} \in G$ which implies that $(k+1) \in I_{l+1} \setminus \{i_{l+1}\}$. Then also $k \in I_{l+1}$ which contradicts the assumption that $k = j_l < i_{l+1}$.
Moreover
$$ \|(z, t) -  (z_k, t_k)\|\leq \hat d(A, \Gamma_x) < \varepsilon$$
and from the construction $B$ we have that $z_k = z'_k$ and $t_k = t'_k$ for $k = j_l$, so
$$
\|(z_k, t_k) -  (z'_k, t'_k)\| = 0.
$$
Hence 
$$\|(z', t') -  (z'_k, t'_k)\| < 2\varepsilon.$$
We also have that
$$\|(z', t') -  (z'_{k+1}, t'_{k+1})\| $$
$$\leq \|(z', t') -  (z, t)\| + \|(z, t) -  (z_{k+1}, t_{k+1})\| + \|(z_{k+1}, t_{k+1}) -  (z'_{k+1}, t'_{k+1})\| < 3\varepsilon.$$
One shows that $\|(z', t') -  (z, t)\| < \varepsilon$ and $\|(z, t) -  (z_{k+1}, t_{k+1})\| < \varepsilon$ in the similar way as above. Assumption $k = j_l < i_{l+1}$ implies that $(k+1) \in I_{l+1} \setminus \{i_{l+1}\}$ so $t_{k+1} \notin G$ and from the construction of set $B$ we have that $t'_{k+1} = t_{k+1}$. Therefore
$$\|(z_{k+1}, t_{k+1}) -  (z'_{k+1}, t'_{k+1})\| = \|f(x)(t) - x(t)\| \leq \| x(\theta_x(t)) - x(\theta_x(t) - ) \| < \varepsilon.$$

$(3^\circ)$ Let $k, k+1 \notin \bigcup_{1 \leq l \leq K} I_l$. Then
$$\|(z', t') -  (z, t)\| = \|f(x)(t) - x(t)\| \leq \| x(\theta_x(t)) - x(\theta_x(t) - ) \| < \varepsilon,$$
$$ \|(z, t) -  (z_k, t_k)\| \leq \hat d(A, \Gamma_x) < \varepsilon$$
and 
$$\|(z_k, t_k) -  (z'_k, t'_k)\| = \|f(x)(t_k) - x(t_k)\| \leq \| x(\theta_x(t_k)) - x(\theta_x(t_k) - ) \| < \varepsilon,$$
so we have that
$$\|(z', t') -  (z'_k, t'_k)\| 
\leq \|(z', t') -  (z, t)\| + \|(z, t) -  (z_k, t_k)\| + \|(z_k, t_k) -  (z'_k, t'_k)\| < 3\varepsilon.$$
The same argument works for $k+1$ thus we also have that
$$\|(z', t') -  (z'_{k+1}, t'_{k+1})\| $$
$$\leq \|(z', t') -  (z, t)\| + \|(z, t) -  (z_{k+1}, t_{k+1})\| + \|(z_{k+1}, t_{k+1}) -  (z'_{k+1}, t'_{k+1})\| < 3\varepsilon.$$

$(4^\circ)$ Let $k+1 = i_l$ for some $l \in \{1, \ldots, K\}$. Observe that if $i_l = j_{l-1}$ then $k, k+1 \in I_{l-1}$ which is treated in $(1^\circ)$. Therefore we assume that $i_l > j_{l-1}$. Then
$$\|(z', t') -  (z, t)\| = \|f(x)(t) - x(t)\| \leq \| x(\theta_x(t)) - x(\theta_x(t) - ) \| < \varepsilon,$$
$$ \|(z, t) -  (z_k, t_k)\| \leq \hat d(A, \Gamma_x) < \varepsilon$$
and 
$$\|(z_k, t_k) -  (z'_k, t'_k)\| = \|f(x)(t_k) - x(t_k)\| \leq \| x(\theta_x(t_k)) - x(\theta_x(t_k) - ) \| < \varepsilon,$$
so we have that
$$\|(z', t') -  (z'_k, t'_k)\| 
\leq \|(z', t') -  (z, t)\| + \|(z, t) -  (z_k, t_k)\| + \|(z_k, t_k) -  (z'_k, t'_k)\| < 3\varepsilon.$$
Similar argumentation and the fact that $(z'_{k+1}, t'_{k+1}) = (z_{k+1}, t_{k+1})$ (see construction of $B$) yield estimate
$$\|(z', t') -  (z'_{k+1}, t'_{k+1})\|$$
$$ \leq \|(z', t') -  (z, t)\| + \|(z, t) -  (z_{k+1}, t_{k+1})\| + \|(z_{k+1}, t_{k+1}) -  (z'_{k+1}, t'_{k+1})\| < 2\varepsilon.$$

As a conclusion of cases $(1^\circ) - (4^\circ)$ we otain that 
$$\hat d(B, \Gamma_y) < 3\varepsilon.$$

The construction of graphs $\Gamma_{y_n}$ and sets $B_n$ is analogous to the construction of graph $\Gamma_y$ and set $B$. In the similar way as above one can show that for all $n \geq n_0$ there holds inequality 
$$\hat d(B_n, \Gamma_{y_n}) < 3\varepsilon.$$

Finally we show that for all $n \geq n_0$
$$d^*(B, B_n) = \max_{1 \leq i \leq m} \{\|(z'_i, t'_i) - (z'_{n,i}, t'_{n,i})\| \ : \ (z'_i, t'_i) \in B, (z'_{n,i}, t'_{n,i})\in B_n \} < 3\varepsilon.$$
First assume that $k \in I_l$ for some $l \in \{1, \ldots, K\}$. Then
$$\| z'_k - z'_{n,k} \| 
= \left\|z_k + \(\frac{z_{j_l} - z_{i_l}}{j_l - i_l}\)\(k - i_l\)  - z_{n,k} - \(\frac{z_{n,j_l} - z_{n,i_l}}{j_l - i_l}\)\(k - i_l\) \right\|$$
$$ = \left\| \(z_{i_l} - z_{n,i_l}\) +  \(z_{j_l} - z_{n,j_l}\)\frac{k - i_l}{j_l - i_l}  + \(z_{n,i_l} - z_{i_l}\)\frac{k - i_l}{j_l - i_l}\right\|$$
$$ \leq 2\|z_{i_l} - z_{n,i_l}\| + \| z_{j_l} - z_{n,j_l} \| = 2\|x_n(\tau_n^{(l)}-) - x(\tau^{(l)}-)\| + \|x_n(\tau_n^{(l)}) - x(\tau^{(l)})\| < 3 \varepsilon ,$$
where the last inequality follows from properties of $G_n$. Similarily
$$| t'_k - t'_{n,k} | 
= \left|t_k + \(\frac{t_{j_l} - t_{i_l}}{j_l - i_l}\)\(k - i_l\)  - t_{n,k} - \(\frac{t_{n,j_l} - t_{n,i_l}}{j_l - i_l}\)\(k - i_l\) \right|$$
$$ = \left| \(t_{i_l} - t_{n,i_l}\) +  \(t_{j_l} - t_{n,j_l}\)\frac{k - i_l}{j_l - i_l}  + \(t_{n,i_l} - t_{i_l}\)\frac{k - i_l}{j_l - i_l}\right|$$
$$ \leq 2|t_{i_l} - t_{n,i_l}| + |t_{j_l} - t_{n,j_l}| $$
$$= 2\left|\sup\{s < \tau_n^{(l)} : x_n(s) \neq x_n(\tau_n^{(l)}-)\} - \sup\{s < \tau^{(l)} : x(s) \neq x(\tau^{(l)}-)\} \right| + |\tau_n^{(l)} - \tau^{(l)}| < 3 \varepsilon$$

Now assume that $k \notin \bigcup_{1 \leq l \leq K} I_l$. Then 
$$ \|(z'_k, t'_k) - (z'_{n,k}, t'_{n,k})\| \leq \|(z'_k, t'_k) - (z_k, t_k)\| + \|(z_k, t_k) - (z_{n,k}, t_{n,k})\| + \|(z_{n,k}, t_{n,k}) - (z'_{n,k}, t'_{n,k})\|.$$
Observe that 
$$\|(z'_k, t'_k) - (z_k, t_k)\| = \|z'_k - z_k\| = \|f(x)(t_k) - x(t_k)\| < \varepsilon $$
and 
$$\|(z'_{n,k}, t'_{n,k}) - (z_{n,k}, t_{n,k})\| = \|z'_{n,k} - z_{n,k}\| = \|f(x_n)(t_{n,k}) - x_n(t_{n,k})\| < \varepsilon.$$
Moreover 
$$
\|(z_k, t_k) - (z_{n,k}, t_{n,k})\| \leq d^*(A, A_n) < \varepsilon.
$$
Then 
$$\|(z'_k, t'_k) - (z'_{n,k}, t'_{n,k})\| < 3 \varepsilon.$$

Thus we have shown that 
$$d^*(B, B_n) < 3 \varepsilon.$$
\ck

\subsection{Continuity of mapping $f$ in the $\JJ_1$ topology.}

\begin{theorem} \label{conv_f_J1}
Let $x, x_n \in \DD([0, \infty), \RR^d)$ fulfil assumptions of Theorem \ref{conv_f}. Assume additionally that $x$ satisfies the condition 
$$ (A) \qquad \forall \ t \geq 0 \quad \exists \ V_t \in [0, \infty) \quad \|x(\theta_x(t)) - x(\theta_x(t)-)\| \leq V_t (\theta_x(t) - \eta_x(t)). $$
Then we also have the convergence $f(x_n) \to f(x)$ in the Skorohod $\JJ_1$ topology.
\end{theorem}

\noindent \textbf{Remark 4.} Condition $(A)$ implies that $x$ has only discontinuities of ``stair'' type. Indeed, if $t$ were the discontinuity point of $x$ which is not the right end of some ``stair'', then by Lemma \ref{lem_eta_theta} we would have that $t = \eta_x(t) = \theta_x(t)$. However, in this case by condition $(A)$ we get
$$0 \leq \|x(\theta_x(t)) - x(\theta_x(t)-)\| = \|x(t) - x(t-)\| \leq V_t (\theta_x(t) - \eta_x(t)) = 0$$
which contradicts that $t \in disc(x)$.

\noindent \textbf{Proof of Theorem \ref{conv_f_J1}.} 
Let $y \equiv f(x)$, $y_n \equiv f(x_n)$. We fix $T \notin disc(y)$ and arbitrary $\varepsilon > 0$. Let integer $K = K(\varepsilon)$, sets $G = G(\varepsilon)$ and $G_n = G_n(\varepsilon)$, integer $m$ and ordered subsets $A \subset \Gamma(x)$, $A_n \subset \Gamma(x_n)$, $B \subset \Gamma(y)$, $B_n \subset \Gamma(y_n)$ be as in proof of Theorem \ref{conv_f}.

As noted in Remark 4. function $x$ has only discontinuities of ``stair'' type. Then the set $B$ is constructed in such way that
$$
\forall \ 1 \leq l \leq K \qquad t'_{i_l} < t'_{i_l + 1} < \ldots < t'_{j_l}. 
$$
Moreover points $(z_k, t_k) \in A$, $k \in \{0, \ldots , m\} \setminus \bigcup_{1 \leq l \leq K} I_l$, may be chosen sach that
$$
\forall \ k_1, k_2 \in \{0, \ldots , m\} \setminus \bigcup_{1 \leq l \leq K} I_l \quad k_1 < k_2 \ \Rightarrow t_{k_1} < t_{k_2}.
$$
Recall that for $k \in \{0, \ldots , m\} \setminus \bigcup_{1 \leq l \leq K} I_l$ points $(z'_k, t'_k) \in B$, we have that $t'_k = t_k$. Thus if for $x$ condition $(A)$ is satisfied then one may choose ordered subset $B \subset \Gamma(y)$ such that
$$
0 = t'_0 < t'_1 < \ldots < t'_m = T.
$$ 
Next observe that assumptions $(i) - (iii)$ together with condition $(A)$ imply that there exists $n_0$ such that for all $n \geq n_0$ all discontinuities of $x_n$ with jump bigger than $\varepsilon$ are of the ``stair'' type. Note that $n_0$ is the same as in proof of Theorem \ref{conv_f}. Hence, applying similar argumentation as above, we may show that for all $n \geq n_0$ ordered subsets $B_n \subset \Gamma(y_n)$  may be chosen in such a way that
$$
0 = t'_{n,0} < t'_{n,1} < \ldots < t'_{n,m} = T.
$$ 
Let us define mappings $\lambda_n \ : [0, T] \mapsto [0, T]$ in the following way
$$
\forall \ s \in [t'_{k-1}, t'_k] \quad \lambda_n(s) \stackrel{def}{=} t'_{n, k-1} + \frac{t'_{n, k} - t'_{n, k-1}}{t'_k - t'_{k-1}}(s - t'_{k-1}),  \quad  k \in \{1, 2, \ldots, m\}.
$$
It is clear that each $\lambda_n$ is strictly increasing and continuous. Moreover
$$
\|\lambda_n - e \| = \max_{0 \leq k \leq m} \|\lambda_n(t'_k) - t'_k\| = \max_{0 \leq k \leq m} \|t'_{n, k} - t'_k\| \leq d^*(B, B_n),
$$
where $e$ denotes identity function. Then, using the estimate of $d^*(B, B_n)$ from the proof of Theorem \ref{conv_f} we get that
$$
\forall \ n \geq n_0 \quad \|\lambda_n - e \| \leq d^*(B, B_n) \leq 3 \varepsilon.
$$

Now we show that $\|y_n \circ \lambda_n - y\| < 9 \varepsilon$ for all $n \geq n_0$. Obviously
$$ 
\|y_n \circ \lambda_n - y\| = \max\left\{ \max_{1 \leq k \leq m} \ \sup_{t \in [t'_{k-1}, t'_k)}\|y_n(\lambda_n(t)) - y(t)\| \ , \ \ \|y_n(\lambda_n(T)) - y(T)\| \right\}
$$
Observe that for any $1 \leq k \leq m$
$$\sup_{t \in [t'_{k-1}, t'_k)}\|y_n(\lambda_n(t)) - y(t)\|$$
$$ = \sup_{t \in [t'_{k-1}, t'_k)} \Big( \|y_n(\lambda_n(t)) - y_n(\lambda_n(t'_{k-1}))\| + \|y_n(\lambda_n(t'_{k-1})) - y(t'_{k-1})\| + \|y(t'_{k-1}) - y(t)\| \Big)$$
$$ \leq \sup_{t \in [t'_{k-1}, t'_k)} \left\|(y_n(\lambda_n(t)), \lambda_n(t)) - (y_n(\lambda_n(t'_{k-1})), \lambda_n(t'_{k-1})) \right\|$$
$$ \qquad + \ \|(y_n(\lambda_n(t'_{k-1})), \lambda_n(t'_{k-1})) - (y(t'_{k-1}), t'_{k-1})\|$$
$$ \qquad + \ \sup_{t \in [t'_{k-1}, t'_k)} \left\|(y(t), t) - (y(t'_{k-1}), t'_{k-1}) \right\|$$
Each mapping $\lambda_n$ is continuous and strictly increasing, thus it transforms interval $[t'_{k-1}, t'_k)$ into $[\lambda_n(t'_{k-1}), \lambda_n(t'_k)) = [t'_{n, k-1}, t'_{n, k})$. Therefore for all $n \geq n_0$ we have that
$$\sup_{t \in [t'_{k-1}, t'_k)} \left\|(y_n(\lambda_n(t)), \lambda_n(t)) - (y_n(\lambda_n(t'_{k-1})), \lambda_n(t'_{k-1})) \right\|$$
$$ = \sup_{t' \in [t'_{n, k-1}, t'_{n, k})} \left\|(y_n(t'), t') - (y_n(t'_{n, k-1}), t'_{n, k-1}) \right\|$$
$$ = \sup_{(z'_{n, k-1}, t'_{n, k-1}) \leq (z', t') < (z'_{n, k}, t'_{n, k}), \ (z', t') \in \Gamma_{y_n}} \left\|(z', t') - (z'_{n, k-1}, t'_{n, k-1}) \right\|$$
$$ \qquad \leq \hat{d}(B_n, \Gamma_{y_n}) < 3 \varepsilon.$$

Similarly, for all $n \geq n_0$  
$$\sup_{t' \in [t'_{n, k-1}, t'_{n, k})} \left\|(y(t'), t') - (y(t'_{k-1}), t'_{k-1}) \right\|$$
$$ = \sup_{(z'_{k-1}, t'_{k-1}) \leq (z', t') < (z'_k, t'_k), \ (z', t') \in \Gamma_{y}} \left\|(z', t') - (z'_{k-1}, t'_{k-1}) \right\|$$
$$ \qquad \leq \hat{d}(B, \Gamma_y) < 3 \varepsilon.$$

Next observe that
$$
\|(y_n(\lambda_n(t'_{k-1})), \lambda_n(t'_{k-1})) - (y(t'_{k-1}), t'_{k-1})\| = \|(z'_{n, k-1}, t'_{n, k-1}) - (z'_{k-1}, t'_{k-1})\| \ , 
$$
where $(z'_{n, k-1}, t'_{n, k-1}) \in B_n$ and$(z'_{k-1}, t'_{k-1}) \in B$. Then, as shown in the proof of Theorem \ref{conv_f}, for $n \geq n_0$ we have that
$$
\|(y_n(\lambda_n(t'_{k-1})), \lambda_n(t'_{k-1})) - (y(t'_{k-1}), t'_{k-1})\| \leq d^*(B, B_n) < 3 \varepsilon.
$$

Finally, for $n \geq n_0$,
$$
|y_n(\lambda_n(T)) - y(T)| = |y_n(T) - y(T)| \leq \|(z'_{n, m}, t'_{n, m}) - (z'_m, t'_m)\| \leq \hat{d}(B, \Gamma_y) < 3 \varepsilon.
$$

Therefore we have shown that for arbitrary small $\varepsilon > 0$ there exist $n_0$ such that for all $n \geq n_0$ we have that
$$
d_{\JJ_1}(y_n, y) \leq \max\{\|y_n \circ \lambda_n - y\|, \|\lambda_n - e\|\} < 9 \varepsilon, 
$$
where $d_{\JJ_1}(\cdot, \cdot)$ is the Skorohod metric generating the $\JJ_1$ topology. This completes the proof. \ck

\section{Functional limit theorems for sequence of CTRW with continuous paths.}
Now let us consider the triangular array of random vectors $\{(Y_{n,k}, J_{n,k})\}$, where $Y_{n,k}$ are the random vectors in $\RR^d$ and $J_{n,k}$ are the positive random variables. For this array we define
\begin{eqnarray}
\(S_n(t), T_n(t) \) &\stackrel{def}{=}& \(\sum_{i=1}^{[nt]} Y_{n,i}, \sum_{i=1}^{[nt]} J_{n,i} \), \quad t\geq 0, \nonumber \\
N_n(t) &\stackrel{def}{=}& \max\{k \geq 0 : T_n(k/n) \leq t \}, \quad t\geq 0, \label{def_Nn}
\end{eqnarray}
and the sequence of CTRW processes 
$$ R_n(t) \stackrel{def}{=} S_n(N_n(t)/n) = \sum_{i=1}^{N_n(t)} Y_{n,i}. $$
In \cite{StrakaHenry} it was shown that processes $R_n$ may be written in the form
$$ R_n = \Phi(S_n, T_n) $$
where the mapping $\Phi: \DD([0,\infty), \mathbb{R}^d \times [0, \infty)) \mapsto \DD([0,\infty), \mathbb{R}^d)$ is given by the formula
$$
\Phi(x, y) = (x^- \circ (y^{-1})^-)^+ .
$$

\begin{proposition}[see Th. 3.6 \cite{StrakaHenry}] \label{prop_conv_CTRW}
Assume that sequence $(S_n, T_n)$ converges weakly in the $\JJ_1$ topology to the process $(A, D)$, where $D$ has almost surely strictly increasing trajectories. Then 
$$
R_n \Rightarrow R = \Phi(A, D)
$$
in the $\JJ_1$ topology.
\end{proposition}

The array $\{(Y_{n,k}, J_{n,k})\}$ generates also the sequence of CPCTRW
$$
\bar R_n(t) \stackrel{df}{=} S_n(N_n(t)/n) + \frac{S_n(N_n(t)/n+1/n) - S_n(N_n(t)/n)}{T_n(N_n(t)/n+1/n) - T_n(N_n(t)/n)}\(t - T_n(N_n(t)/n)\) \quad t\geq 0.
$$
Observe that $\eta_{R_n}(t) = \sup\{s < t: R_n(s) \neq R_n(t)\}$ and $\theta_{R_n}(t) = \inf\{s \geq t: R_n(s) \neq R_n(t-)\}$ are the moments of the last jump of $R_n$ before time $t$ and the first jump after time $t$, respectively. From the definition (\ref{def_Nn}) of the counting process $N_n$ it then follows that
$$
\eta_{R_n}(t) = T_n(N_n(t)/n) \quad \mbox{and} \quad \theta_{R_n}(t) = T_n(N_n(t)/n + 1/n).
$$
Then we also have that
$$
R_n(\eta_{R_n}(t)) = S_n\(N_n(\eta_{R_n}(t))/n \) = S_n\(N_n(T_n(N_n(t)/n))/n \).
$$
On the other hand, from definition (\ref{def_Nn}) we get
$$
N_n(T_n(N_n(t)/n)) = \max\{k \geq 0 : T_n(k/n) \leq T_n(N_n(t)/n)\} = N_n(t).
$$
Putting these two things together we obtain
$$
R_n(\eta_{R_n}(t)) = S_n\(N_n(t)/n\).
$$
Similarily
$$
R_n(\theta_{R_n}(t)) = S_n\(N_n(\theta_{R_n}(t))/n \) = S_n\(N_n(T_n(N_n(t)/n) + 1/n)/n \).
$$
and
$$
N_n(T_n(N_n(t)/n) + 1/n) = \max\{k \geq 0 : T_n(k/n) \leq T_n(N_n(t)/n + 1/n)\} = N_n(t) + 1
$$
thus we get
$$
R_n(\theta_{R_n}(t)) = S_n\(N_n(t)/n + 1/n\).
$$

Hence we can write
$$
\bar R_n(t) = R_n(\eta_{R_n}(t)) + \frac{R_n(\theta_{R_n}(t)) - R_n(\eta_{R_n}(t))}{\theta_{R_n}(t) - \eta_{R_n}(t)}\(t - \eta_{R_n}(t)\) = f(R_n)(t),
$$ 
where $f$ is given by (\ref{def_f}). Below we state and prove the functional limit theorem for sequence $\bar R_n$.

\begin{theorem} \label{th_conv_CPCTRW}
Assume that $(S_n, T_n) \Rightarrow (A, D)$ in the space $\DD\([0, \infty), \RR^d \times [0, \infty) \)$ equipped with $\JJ_1$ topology, where almost surely the trajectories of the process $A$ are not constant on any interval $(a, b) \in [0, \infty)$ and $D$ has almost surely strictly increasing realizations. Then:

\begin{itemize}
	\item[(a)] 
	\begin{equation} \label{conv_bar_Rn}
	\bar R_n \Rightarrow \bar R \stackrel{def}{=} f(R) = f(\Phi(A, D)) 
	\end{equation}
in the space $\DD\([0, \infty), \RR^d \)$ equipped with the $S\MM_1$ topology.
	\item[(b)] If additionally 
	$$P(disc(A) \subset disc(D)) = 1$$
	then convergence (\ref{conv_bar_Rn}) holds in the Skorohod $\JJ_1$ topology.
\end{itemize}

\end{theorem}

In the proof of Theorem \ref{th_conv_CPCTRW} we will make use of the following lemma.
\begin{lemma} \label{Phi_equiv}
For any function $y \in \DD\([0, \infty), [0, \infty) \)$ which is nondecreasing and unbounded from above the following equality holds
$$
(y \circ y^{-1})^{-1} = \left(y^- \circ (y^{-1})^-\right)^+ \equiv \Phi(y, y).
$$
\end{lemma}
For the proof of lemma see Appendix.\\

\noindent \textbf{Proof of Theorem \ref{th_conv_CPCTRW}.} Part $(a)$. By Skorohod representation (see eg. Th. 3.2.2 \cite{Whitt}) there exists some probability space $\Omega'$ and sequence of processes $(S_n', T_n')$ and process $(A', D')$ defined on this space and such that $(S_n', T_n') \stackrel{d}{=} (S_n, T_n)$, $(A', D') \stackrel{d}{=} (A, D)$ and $(S_n', T_n') \to (A' D')$ almost surely in the $\JJ_1$ topology. Then also $R_n' \stackrel{def}{=} \Phi(S_n', T_n') \stackrel{d}{=} R_n$, $R' \stackrel{def}{=} \Phi(A', D') \stackrel{d}{=} R$ and $R_n' \to R'$ almost surely in the $\JJ_1$ topology.

Now we check that processes $R_n'$, $R'$ satisfy assumptions of Theorem \ref{conv_f}. Assume that $R'$ is constant on $(\eta_{R'}(t), \theta_{R'}(t))$ for some $t > 0$ and $\theta_{R'}(t) \in disc(R')$. Note that $A'$ is not constant on any interval $(a, b) \in [0, \infty)$ so it is easy to see that $R' = \Phi(A', D')$ is constant exactly when $\Phi(D', D')$ is constant. By Lemma \ref{Phi_equiv} $$\Phi(D', D') = \(D' \circ D'^{-1}\)^{-1}$$
and one easily checks (see sect. 3 of \cite{FCCTRW}) that the left side of this equality is constant 
in some neighbourhood of $t$ if and only if
$$ 
t \in [D'(\tau-), D'(\tau)) \quad \mbox{for some } \ \tau \in disc(D').
$$ 
Then

\begin{equation} \label{eq_3a}
\eta_{R'}(t) = D'(\tau-) \ , \quad \theta_{R'}(t) = D'(\tau).
\end{equation}

Convergence $(S_n', T_n') \to (A', D')$ in the $\JJ_1$ topology implies that there exist sequence $\tau_n \to \tau$ with $\tau \in disc((A', D'))$ and $\tau_n \in disc((S_n',T_n'))$ for all $n \geq n_0$ for some $n_0$, for which sequence the following convergences hold
\begin{equation} \label{eq_3b}
\(S_n'(\tau_n - ), T_n'(\tau_n - )\) \to \(A'(\tau-), D'(\tau-)\) \ , \ \ \(S_n'(\tau_n), T_n'(\tau_n)\) \to \(A'(\tau), D'(\tau)\)
\end{equation}
(see e.g. argument on page 79 of \cite{Whitt}).

As $T_n'(\tau_n - ) \to D'(\tau-)$, $T_n'(\tau_n) \to D'(\tau)$ and $\tau \in disc(D')$, we may then choose sequence $t_n \to t$ such that $\forall n \geq n_0$ $t_n \in \(T_n'(\tau_n - ), T_n'(\tau_n)\)$. From the definition of $R_n'$ it easily follows that they are constant on intervals $[T_n'(\tau_n - ), T_n'(\tau_n))$ and we get that
\begin{equation} \label{eq_3c}
\eta_{R'_n}(t_n) = T_n'(\tau_n-) \ , \quad \theta_{R_n'}(t_n) = T_n'(\tau).
\end{equation}
Form equations (\ref{eq_3a}) - (\ref{eq_3c}) we obtain convergences
$$
\eta_{R_n'}(t_n) \to \eta_{R'}(t) \ , \quad \mbox{and} \quad \theta_{R_n'}(t_n) \to \theta_{R'}(t),
$$
so the processes $R_n'$, $R'$ satisfy assumptions $(i)$ and $(ii)$ of Theorem \ref{conv_f}. Now we check that assumption $(iii)$ is also satisfied.

First notice that by definition for any $n \geq 1$ both processes $S_n'$ and $T_n'$ are constant on the intervals of the form $[i/n, (i+1)/n)$, $i \in \{0, 1, 2, \ldots\}$ and so $disc(S_n') = disc(T_n')$. For arbitrary fixed $n$, by definition of $R_n'$ and (\ref{eq_3c}) we may write
\begin{eqnarray*}
R_n'\(\theta_{R_n'}(t_n)\) &=& \(S_n'^- \circ (T_n'^{-1})^-\)^+(T_n'(\tau_n)) \\
&=& \lim_{h \searrow 0} \lim_{u \searrow 0} \lim_{v \searrow 0} S_n'\(T_n'^{-1}(T_n'(\tau_n) + h - v) - u \)
\end{eqnarray*}
where the limits are resolved starting from the most inner one. 
Since $T_n'$ is a step function constant on intervals of length $1/n$ it follows that 
$$
T_n'^{-1}(T_n'(\tau_n)) = \inf\{s>0 : T_n'(s) > T_n'(\tau_n)\} = \tau_n + 1/n \in disc(T_n'),
$$
i.e. it is the first moment of jump of $T_n'$ after time $\tau_n$. Moreover $T_n'^{-1}$ is constant on interval $[T_n'(\tau_n), T_n'(\tau_n + 1/n)) = [T_n'(\tau_n + 1/n -), T_n'(\tau_n + 1/n))$ so for sufficiently small $h$ and $v < h$ we get that
$$
T_n'^{-1}(T_n'(\tau_n) + h - v) = \tau_n + 1/n
$$
and then
$$
R_n'\(\theta_{R_n'}(t_n)\) = \lim_{u \searrow 0} S_n'\(\tau_n + 1/n - u \).
$$
However $S_n'$ is also a step function and is constant on the interval $[\tau_n, \tau_n + 1/n)$ so for $u < 1/n$ we have that $S_n'\(\tau_n + 1/n - u \) = S_n'\(\tau_n\)$ and
\begin{equation} \label{eq_3d}
R_n'\(\theta_{R_n'}(t_n)\) = S_n'\(\tau_n\).
\end{equation}
Analogous argumentation leads to
\begin{eqnarray}
R_n'\(\theta_{R_n'}(t_n) - \) &=& \(S_n'^- \circ (T_n'^{-1})^-\)(T_n'(\tau_n)) \nonumber \\ 
&=&  \lim_{u \searrow 0} \lim_{v \searrow 0} S_n'\(T_n'^{-1}(T_n'(\tau_n) - v) - u \) \nonumber \\ 
&=&  \lim_{u \searrow 0} S_n'\(\tau_n - u \) = S_n'\(\tau_n - \). \label{eq_3e}
\end{eqnarray}
To check assumption $(iii)$ we also need to find $R'(\theta_{R'}(t))$ and $R'(\theta_{R'}(t)-)$. Recall that we assume $D'$ to be strictly increasing, hence $D'^{-1}$ is continuous and $(D'^{-1})^- = D'^{-1}$. By definition of $R'$ and equation (\ref{eq_3a}) we have
$$
R'(\theta_{R'}(t)) = \(A'^- \circ D'^{-1}\)^+(D'(\tau)) = \lim_{h \searrow 0} \lim_{u \searrow 0} A'\(D'^{-1}(D'(\tau) +h) - u \).
$$
Since $D'$ is strictly increasing and right-continuous it follows that
$$
D'^{-1}(D'(\tau)) = \inf\{s > 0 : D'(s) > D'(\tau)\} = \tau
$$
and there exists $\varepsilon(h) \searrow 0$ as $h \searrow 0$ such that $D'^{-1}(D'(\tau) + h) = \tau + \varepsilon(h)$. Then right-continuity of $A'$ yields
\begin{equation} \label{eq_3f}
R'(\theta_{R'}(t)) = \lim_{h \searrow 0} \lim_{u \searrow 0} A'\(\tau + \varepsilon(h) - u \) = \lim_{h \searrow 0} 
 A'\(\tau + \varepsilon(h)\) = A'(\tau).
\end{equation}
We also have that
\begin{eqnarray} 
R'(\theta_{R'}(t)-) &=& \lim_{u \searrow 0} \(A'^- \circ D'^{-1}\)(D'(\tau)) = \lim_{u \searrow 0} A'\(D'^{-1}(D'(\tau)) - u \) \nonumber \\
 &=& \lim_{u \searrow 0} A'\(\tau - u \) = A'(\tau-). \label{eq_3g}
\end{eqnarray}
Combining convergence (\ref{eq_3b}) with equations (\ref{eq_3d})-(\ref{eq_3g}) we get 
$$
\|R_n'(\theta_{R_n'}(t)) - R_n'(\theta_{R_n'}(t_n)-)\| \to \|R'(\theta_{R'}(t)) - R'(\theta_{R'}(t)-)\| \quad \mbox{as} \ n \to \infty 
$$
and the assumption $(iii)$ of Theorem \ref{conv_f} is satisfied. Then 
$$
\bar R_n' \stackrel{def}{=} f(R_n') \to \bar R' \stackrel{def}{=} f(R') \quad \mbox{almost surely in the } S\MM_1 \mbox{ topology}
$$
and as $\bar R_n' \stackrel{d}{=} \bar R_n$ and $\bar R' \stackrel{d}{=} \bar R$ there also follows weak convergence
$$
f(R_n) \Rightarrow f(R) \quad \mbox{in the } S\MM_1 \mbox{ topology}
$$
which completes the proof of part $(a)$. \\

The proof of patr $(b)$ of the theorem follows in the same way as proof part $(a)$ with Theorem \ref{conv_f_J1} used in place of Theorem \ref{conv_f}. Thus we only need to check if paths of $R'$ satisfy condition $(A)$ almost surely. 

Recall, that paths of $A'$ are nowhere constant almost surely (see assumptions). Then $R' = \Phi(A', D')$ is constant exactly on the same intervals as $\Phi(D', D')$, that is on intervals $[D'(\tau -), D'(\tau))$, $\tau \in disc(D')$. Moreover, assumption that 
$disc(A') \subset disc(D')$ implies that $$disc(R') \subset disc(\Phi(D', D')) = \{D'(\tau), \ \tau \in disc(D') \}$$  
(see sect. 3 of \cite{FCCTRW}).

If $t \in [D'(\tau -), D'(\tau))$ for some $\tau \in disc(D')$ then by equalities (\ref{eq_3a}),(\ref{eq_3f}) and (\ref{eq_3g}) we have that 
$$
\|R'(\theta_{R'}(t)) - R'(\theta_{R'}(t)-)\| = \|A'(\tau) - A'(\tau-)\| 
$$
and 
$$
|\theta_{R'}(t) - \eta_{R'}(t)| = |D'(\tau) - D'(\tau-)|
$$ 
and one easily finds $V_t \in [0, \infty)$ such that condition $(A)$ is satisfied.

If $t \notin [D'(\tau -), D'(\tau))$ for any $\tau \in disc(D')$, then $R'$ is not constant in the neighbourhood of $t$ and by Lemma \ref{lem_eta_theta} $\eta_{R'}(t) = \theta_{R'}(t) = t$. Since $t \notin disc(R') \subset  \{D'(\tau), \ \tau \in disc(D') \}$ it follows that
$$
\|R'(\theta_{R'}(t)) - R'(\theta_{R'}(t)-)\| = \|R'(t) - R'(t-)\| = 0.
$$
Then condition $(A)$ holds with any $V_t \in [0, \infty)$. This completes the proof of part $(b)$. \ck

\section*{Appendix}
\textbf{Proof of Lemma \ref{Phi_equiv}.}
We split the halfline $[0, \infty)$ into two sets $B_1$ i $B_2$ where
$$B_1 \stackrel{def}{=} \bigcup_{\tau_y \in disc(y)}[y(\tau_y -), y(\tau_y)) \quad \mbox{and} \quad 
B_2 \stackrel{def}{=} [0, \infty) \setminus B_1.$$
First let us assume that $t \in B_1$. Then $t \in [y(\tau_y^*-), y(\tau_y^*))$ for some $\tau_y^* \in Disc(y)$ and there exists $h_0 > 0$ such that $\forall \ h \in(0, h_0)$
$$t+h \in (y(\tau_y^*-), y(\tau_y^*)).$$
Then
$$
\left(y^{-1}\right)^-(t+h) = \lim_{a \searrow 0} \left(y^{-1}\right)((t+h)-a) = \lim_{a \searrow 0} \inf\{s>0 : y(s) > (t+h)-a\} $$
$$ = \inf\{s>0 : y(s) > y(\tau_y^*-)\} = \tau_y^* \ 
$$
and it follows that
$$
\lim_{h \searrow 0} \left(y^- \circ (y^{-1})^-\right)(t+h) = y(\tau_y^*-).
$$
Hence we obtain that for all $t \in [y(\tau_y^*-), y(\tau_y^*))$
$$
\left(y^- \circ (y^{-1})^-\right)^+(t) = \lim_{h \searrow 0} \left(y^- \circ (y^{-1})^-\right)(t+h) 
= y(\tau_y^*-).
$$
On the other hand, by formula (13) in \cite{FCCTRW} we get that for any $t \in [y(\tau_y^*-), y(\tau_y^*))$
$$
\left(y \circ y^{-1}\right)^{-1}(t) = y(\tau_y^*-).
$$
Joining these two observations we obtain that
$$
\left(y \circ y^{-1}\right)^{-1}(t) = \left(y^- \circ (y^{-1})^-\right)^+(t) = \Phi(y, y) \quad  \forall \ t \in B_1.
$$

Now assume that $t \in B_2$ and let $u_0 = \sup\{s : y(s) = t \}$. Then $y(u_0) = t$, which is obvious in the case when there is only one $s \geq 0$ such that $y(s) = t$. If $y(s) = t$ for more than one $s$, then $y(s_0 - ) = t$ which together with our assumption on $t$ being in the set $B_2$ implies equality $y(u_0) = t$. Since $y$ is nondecreasing and continuous from the right, we may choose sequence $u_n \searrow u_0$ such that $y(u_n) \searrow y(u_0)$. Then for sequence $\{u_n\}$ we define sequence $h_n \searrow 0$ such that $t + h_n = y(u_n)$ and so we may write
$$
\left(y^- \circ (y^{-1})^-\right)^+(t) = \lim_{h \searrow 0} \left(y^- \circ (y^{-1})^-\right)(t+h) = 
\lim_{n \to \infty} \left(y^- \circ (y^{-1})^-\right)(t+h_n).
$$ 
Note however that
$$
\left(y^{-1}\right)^-(t+h_n) = \lim_{a \searrow 0} \inf\{s>0 : y(s) > (t+h_n)-a\} \in [u_{n+1}, u_n].
$$
Then
$$
\left(y^- \circ (y^{-1})^-\right)(t+h_n) \in [y^-(u_{n+1}), y^-(u_n)]
$$
and, as $n \to \infty$, we get that
$$
t = y(u_0) \leq y^-(u_{n+1}) \leq \left(y^- \circ (y^{-1})^-\right)(t+h_n) \leq y^-(u_n) \leq y(u_n) \to t.
$$
Thus $\left(y^- \circ (y^{-1})^-\right)^+(t) = t$ for $t \in B_2$. It turns out that for such $t$ we also have that $\left(y \circ y^{-1}\right)^{-1}(t) = t$ (see argumentation justifying that formula (14) in \cite{FCCTRW} holds true). Combining these two equalities we get that
$$\left(y \circ y^{-1}\right)^{-1}(t) = \left(y^- \circ (y^{-1})^-\right)^+(t) = \Phi(y, y)(t) \quad  \forall \ t \in B_2$$
and the proof is complete.
\ck

\textbf{Acknowledgments.}\\
The research of M. Magdziarz was partially supported by NCN Maestro grant no. 2012/06/A/ST1/00258.


\end{document}